\documentclass[a4paper,10pt,twoside]{article}
\usepackage{amsmath,amssymb,amsthm,eucal,graphics}
\usepackage[cp1251]{inputenc}
\usepackage[russian]{babel}
\usepackage{graphicx}
\usepackage{amsfonts,amssymb,eucal}
\usepackage{amsbsy}
\usepackage{latexsym}
\usepackage{tabularx}

\makeatletter \gdef\firstpage{1}

\def\frsthdr{ }

\def\firstpageone{0\thepage}
\def\firstpagetwo{00\thepage}
\def\firstpagethree{000\thepage}
\def\firstpagemark{\ifnum\firstpage <10  \firstpageone
\else\ifnum\firstpage<100 \firstpagetwo \else \ifnum\firstpage
<1000 \firstpagethree \else \firstpageone\fi\fi\fi}

\def\footline{\ifnum\thepage=\firstpage \footlineone
\else\footlinetwo\fi}
\def\footlineone{ }
\def\footlinetwo{}

\def\titles{Finite almost simple groups}
\def\authors{ A.\,A.\,Makhnev}

\def\oddhedr{\ifnum\thepage=\firstpage \firsthdr \else \odhdr \fi}
\def\firsthdr{\hspace{\fill} \sl \frsthdr \hspace{\fill}\hbox{}}
\def\odhdr{\hspace{\fill}\sl\rightmark \titles \hspace{\fill}
\rm \thepage}

\def\evnhedr{\ifnum\thepage=\firstpage \firsthdr \else \evhdr \fi}
\def\evhdr{\noindent \rm \thepage\hspace*{\fill} \sl\leftmark
\authors \hspace*{\fill}\hbox{}}

\def\ps@newpstyle{\def\@oddhead{
\hspace{-0.65em} \vbox{\oddhedr\vskip 1mm \hrule width \textwidth}
}
\def\@evenhead{
\hspace{-0.65em} \vbox{\evnhedr\vskip 1mm \hrule width \textwidth}
}\textsc{}
\def\@oddfoot{\footline}
\def\@evenfoot{\@oddfoot}
}

\def\refer
{
\begin{small}
\begin{enumerate}
\leftmargin=0pt \rightmargin=0pt \itemsep=0pt
\parsep=0pt
}

\def\endref{\end{enumerate}\end{small} }

\begin{document}

\setcounter{page}{\firstpage} 

\Russian \sloppy \rm


\newcommand{\NN}{\mathbb{N}}
\newcommand{\ZZ}{\mathbb{Z}}
\newcommand{\RR}{\mathbb{R}}

\newcommand{\inN}{\in\NN}
\newcommand{\inZ}{\in\ZZ}
\newcommand{\inR}{\in\RR}

\newcommand{\alf}{\alpha}
\newcommand{\bet}{\beta}
\newcommand{\eps}{\varepsilon}
\newcommand{\lam}{\lambda}
\newcommand{\ups}{\upsilon}
\newcommand{\sgn}{\rm sign}

\newtheorem{teo}{\bf ╥┼╬╨┼╠└}
\newtheorem{lem}{\bf ╦┼╠╠└}
\newtheorem{hyp}{\bf ╧╨┼─╧╬╦╬╞┼═╚┼}
\newtheorem{prp}{\bf ╧╨┼─╦╬╞┼═╚┼}
\newtheorem{ass}{\bf ╙╥┬┼╨╞─┼═╚┼}
\newtheorem{cor}{\bf ╤╦┼─╤╥┬╚┼}

\newcommand{\pr}{\par\mbox{─╬╩└╟└╥┼╦▄╤╥┬╬.}\ \ }

\renewcommand{\theteo}{\arabic{teo}}
\renewcommand{\theprp}{\arabic{prp}}
\renewcommand{\theass}{\arabic{ass}}
\renewcommand{\thelem}{\arabic{lem}}
\renewcommand{\thecor}{\arabic{cor}}
\renewcommand{\thehyp}{\arabic{hyp}}

\begin{center} MOORE GRAPH WITH PARAMETERS (3250,57,0,1) DOES NOT EXIST\\
(in Russian, English abstract)\\
A.\,A.\,Makhnev$^{1,2}$\\
${^1}$Krasovskii Institute of Mathematics and Mechanics UB RAS, \\
${^2}$Ural Federal University\\
makhnev@imm.uran.ru
\end{center}

\allowdisplaybreaks
\renewcommand{\tablename}{}


{\small \small \emph{Abstract.} If a regular graph of degree $k$ and diameter $d$ has $v$ vertices then
$$v\le 1+k+k(k-1)+\dots+k(k-1)^{d-1}.$$
Graphs with $v=1+k+k(k-1)+\dots+k(k-1)^{d-1}$ are called Moore graphs. Damerell proved that a Moore graph of degree $k\ge 3$
has diameter $2$. If $\Gamma$ is a Moore graph of diameter $2$, then $v=k^2+1$, $\Gamma$ is strongly regular with $\lambda=0$ and $\mu=1$, and one of the following statements holds{\rm:} $k=2$ and $\Gamma$ is the pentagon, $k=3$ and $\Gamma$ is the Petersen graph, $k=7$ and $\Gamma$ is the Hoffman-Singleton graph, or $k=57$. The existence of a Moore graph of degree $57$ was unknown.

Jurishich and Vidali have proved that the existence of a Moore graph of degree $k>3$ is equivalent to the existence of a distance-regular graph with intersection array $\{k-2,k-3,2;1,1,k-3\}$ (in the case $k=57$ we have a distance-regular graph with intersection array $\{55,54,2;1,1,54\}$).
\smallskip

In this paper we prove that a distance-regular graph with intersection array $\{55,54,2;1,1,54\}$ does not exist.
As a corollary, we prove that a Moore graph of degree $57$ does not exist.

}

\section*{┬тхфхэшх}

╠√ ЁрёёьрЄЁштрхь эхюЁшхэЄшЁютрээ√х уЁрЇ√ схч яхЄхы№ ш ъЁрЄэ√ї ЁхсхЁ.
┼ёыш $a,b$ --- тхЁ°шэ√ уЁрЇр $\Gamma$, Єю ўхЁхч $d(a,b)$ юсючэрўрхЄё 
ЁрёёЄю эшх ьхцфє $a$ ш $b$, р ўхЁхч $\Gamma_{i}(a)$ --- яюфуЁрЇ уЁрЇр
$\Gamma$, шэфєЎшЁютрээ√щ ьэюцхёЄтюь тхЁ°шэ, ъюЄюЁ√х эрїюф Єё  эр ЁрёёЄю эшш
$i$ т $\Gamma$ юЄ тхЁ°шэ√ $a$. ╧юфуЁрЇ $\Gamma_1(a)$ эрч√трхЄё  {\it
юъЁхёЄэюёЄ№■ тхЁ°шэ√} $a$ ш юсючэрўрхЄё  ўхЁхч $[a]$, хёыш уЁрЇ $\Gamma$ ЇшъёшЁютрэ.

╬яЁхфхыхэшх фшёЄрэЎшюээю Ёхуєы Ёэюую уЁрЇр ё ьрёёштюь яхЁхёхўхэшщ $\{b_0,...,b_{d-1};\linebreak c_1,...,c_d\}$
ьюцэю эрщЄш т \cite{BCN}.

─ы  уЁрЇр $\Gamma$ фшрьхЄЁр $d$ ш $i\in \{2,...,d\}$ ўхЁхч $\Gamma_i$ юсючэрўшь уЁрЇ
эр Єюь цх ьэюцхёЄтх тхЁ°шэ, т ъюЄюЁюь тхЁ°шэ√ $u,w$ ёьхцэ√, хёыш ЁрёёЄю эшх ьхцфє
эшьш т  $\Gamma$ Ёртэю $i$.

┼ёыш Ёхуєы Ёэ√щ уЁрЇ ёЄхяхэш $k$ ш фшрьхЄЁр $d$ шьххЄ $v$ тхЁ°шэ, Єю т√яюыэ хЄё  эхЁртхэёЄтю:
$$v\le 1+k+k(k-1)+\dots+k(k-1)^{d-1}.$$
├ЁрЇ√, фы  ъюЄюЁ√ї ¤Єю эхёЄЁюуюх эхЁртхэёЄтю яЁхтЁр∙рхЄё  т ЁртхэёЄтю, эрч√тр■Єё  уЁрЇрьш ╠єЁр (1960). ▌Єю юяЁхфхыхэшх
яЁшэрфыхцшЄ ╒юЇьрэє ш ╤шэуыЄюэє \cite{HS}, ъюЄюЁ√х юяшёрыш уЁрЇ√ ╠єЁр фшрьхЄЁр 2 ш 3.  ╧ЁюёЄхщ°шщ яЁшьхЁ уЁрЇр ╠єЁр
фюёЄрты хЄ $(2d+1)$-єуюы№эшъ.

─рьхЁхыы \cite{D} (ёь. Єръцх \cite{BI}) фюърчры, ўЄю уЁрЇ ╠єЁр ёЄхяхэш $k\ge 3$ шьххЄ фшрьхЄЁ 2 . ┬ ¤Єюь ёыєўрх $v=k^2+1$,
уЁрЇ ёшы№эю Ёхуєы Ёхэ ё $\lambda=0$ ш $\mu=1$, р ёЄхяхэ№ $k$ Ёртэр 2 (я Єшєуюы№эшъ), 3 (уЁрЇ ╧хЄхЁёхэр), 7 (уЁрЇ ╒юЇьрэр-╤шэуыЄюэр)
шыш 57. ╤є∙хёЄтютрэшх уЁрЇр ╠єЁр ёЄхяхэш $k=57$ Ёрэхх с√ыю эхшчтхёЄэю.

└°срїхЁ \cite{A} фюърчры, ўЄю уЁрЇр ╠єЁр ёЄхяхэш $57$ эх  ты хЄё  фшёЄрэЎшюээю ЄЁрэчшЄштэ√ь (яю¤Єюьє уЁрЇ ╠єЁр ёЄхяхэш $57$
шэюуфр эрч√тр■Є уЁрЇюь └°срїхЁр). ├. ╒шуьхэ (ёь. \cite[ЄхюЁхьр 3.13]{C}) фюърчры, ўЄю уЁрЇ ╠єЁр ёЄхяхэш $57$ эх  ты хЄё  тхЁ°шээю
ёшььхЄЁшўэ√ь. └.└. ╠рїэхт ш ─.┬. ╧рфєўшї \cite{MP1} ш \cite{MP2} шчєўшыш тючьюцэ√х ртЄюьюЁЇшчь√ уЁрЇр ╠єЁр ёЄхяхэш $57$.

▐Ёш°шў ш ┬шфрыш \cite{JV} чрьхЄшыш, ўЄю ёє∙хёЄтютрэшх уЁрЇр ╠єЁр ёЄхяхэш $k>3$
Ёртэюёшы№эю ёє∙хёЄтютрэш■ фшёЄрэЎшюээю Ёхуєы Ёэюую уЁрЇр ё ьрёёштюь яхЁхёхўхэшщ $\{k-2,k-3,2;1,1,k-3\}$
(т ёыєўрх $k=57$ яюыєўшь уЁрЇ ё ьрёёштюь яхЁхёхўхэшщ $\{55,54,2;1,1,54\}$). └.└. ╠рїэхт ш ─.┬. ╧рфєўшї \cite{MP3} эр°ыш
тючьюцэ√х ртЄюьюЁЇшчь√ фшёЄрэЎшюээю Ёхуєы Ёэюую уЁрЇр ё ьрёёштюь яхЁхёхўхэшщ $\{55,54,2;1,1,54\}$.

┬ ЁрсюЄх шчєўхэ√ ётющёЄтр уЁрЇр ё ьрёёштюь яхЁхёхўхэшщ $\{55,54,2;1,1,54\}$, ё яюью∙№■ ъюЄюЁ√ї яюыєўхэр
\smallskip

\begin{teo}\label{Jgraph} ─шёЄрэЎшюээю Ёхуєы Ёэ√щ уЁрЇ ё ьрёёштюь яхЁхёхўхэшщ $\{55,54,2;1,1,54\}$ эх ёє∙хёЄтєхЄ.
\end{teo}

┼ёЄхёЄтхээ√ь юсЁрчюь яюыєўрхЄё 
\smallskip

\begin{cor}\label{Mgraph} ╤шы№эю Ёхуєы Ёэ√щ уЁрЇ ╠єЁр ё ярЁрьхЄЁрьш $(3250,57,0,1)$ эх ёє∙хёЄтєхЄ.
\end{cor}

▌ЄюЄ Ёхчєы№ЄрЄ ёЄртшЄ чръы■ўшЄхы№эє■ Єюўъє т 60-ыхЄэшї шёёыхфютрэш ї уЁрЇют ╠єЁр.

╬ёэютэ√х ьхЄюф√, шёяюы№чютрээ√х т фюърчрЄхы№ёЄтх яюыєўхээ√ї Ёхчєы№ЄрЄют:

ьхЄюф ЄЁющэ√ї ўшёхы яхЁхёхўхэшщ, яЁшьхэхээ√щ ъ уЁрЇє схч єёыютш  $Q$-яюышэю\-ьшры№\-эюёЄш;

ьхЄюф ёшььхЄЁшчрЎшш ьрёёштр ЄЁющэ√ї ўшёхы яхЁхёхўхэшщ (яЁхфыюцхэ └.└. ╠рїэхт√ь, рэрыюу ёшььхЄЁшчрЎшш ЄхэчюЁют).

\section*{┬ёяюьюурЄхы№э√х Ёхчєы№ЄрЄ√}

═р°ш юсючэрўхэш  ш ЄхЁьшэюыюуш  т юёэютэюь ёЄрэфрЁЄэ√, шї ьюцэю
эрщЄш т \cite{BCN}.

┬ фюърчрЄхы№ёЄтх ЄхюЁхь√ шёяюы№чє■Єё  ЄЁющэ√х ўшёыр яхЁхёхўхэшщ \cite{JV1}.

╧єёЄ№ $\Gamma$ --- фшёЄрэЎшюээю Ёхуєы Ёэ√щ уЁрЇ фшрьхЄЁр $d$. ┼ёыш $u_1,u_2,u_3$ --- тхЁ°шэ√
уЁрЇр $\Gamma$, $r_1,r_2,r_3$ --- эхюЄЁшЎрЄхы№э√х Ўхы√х ўшёыр, эх сюы№°шх $d$.
╫хЁхч $\left\{u_1u_2u_3\atop r_1r_2r_3\right\}$ юсючэрўшь ьэюцхёЄтю тхЁ°шэ $w\in \Gamma$ Єръшї, ўЄю $d(w,u_i)=r_i$,
р ўхЁхч $\left[u_1u_2u_3\atop r_1r_2r_3\right]$ --- ўшёыю тхЁ°шэ т $\left\{u_1u_2u_3\atop r_1r_2r_3\right\}$.
╫шёыр $\left[u_1u_2u_3\atop r_1r_2r_3\right]$ эрч√тр■Єё  ЄЁющэ√ьш ўшёырьш яхЁхёхўхэшщ.
─ы  ЇшъёшЁютрээющ ЄЁющъш тхЁ°шэ $u_1,u_2,u_3$ тьхёЄю $\left[u_1u_2u_3\atop r_1r_2r_3\right]$ сєфхь яшёрЄ№
$[r_1r_2r_3]$. ╩ ёюцрыхэш■, фы  ўшёхы $[r_1r_2r_3]$ эхЄ юс∙шї ЇюЁьєы. ╬фэръю, т \cite{JV1} яЁхфыюцхэ ьхЄюф
т√ўшёыхэш  эхъюЄюЁ√ї ўшёхы $[r_1r_2r_3]$.

╧єёЄ№ $u,v,w$ --- тхЁ°шэ√ уЁрЇр $\Gamma$, $W=d(u,v),U=d(v,w),\\V=d(u,w)$. ╥ръ ъръ шьххЄё  Єюўэю юфэр тхЁ°шэр
$x=u$ Єрър , ўЄю $d(x,u)=0$, Єю ўшёыю $[0jh]$ Ёртэю 0 шыш 1. ╬Єё■фр $[0jh]=\delta_{jW}\delta_{hV}$.
└эрыюушўэю, $[i0h]=\delta_{iW}\delta_{hU}$ ш $[ij0]=\delta_{iU}\delta_{jV}$.

─Ёєуюх ьэюцхёЄтю єЁртэхэшщ ьюцэю яюыєўшЄ№, ЇшъёшЁє  ЁрёёЄю эшх ьхцфє фтєь  тхЁ°шэрьш шч $\{u,v,w\}$
ш ёюёўшЄрт ўшёыю тхЁ°шэ, эрїюф ∙шїё  эр тёхї тючьюцэ√ї ЁрёёЄю эш ї юЄ ЄЁхЄ№хщ:

$$\sum_{l=1}^d[ljh]=p_{jh}^U-[0jh], \sum_{l=1}^d[ilh]=p_{ih}^V-[i0h], \sum_{l=1}^d[ijl]=p_{ij}^W-[ij0]\, (+)$$

╧Ёш ¤Єюь эхъюЄюЁ√х ЄЁющъш шёўхчр■Є. ╧Ёш $|i-j|>W$ шыш $i+j<W$ шьххь $p_{ij}^W=0$, яю¤Єюьє $[ijh]=0$
фы  тёхї $h\in \{0,...,d\}$.

╧юыюцшь $S_{ijh}(u,v,w)=\sum_{r,s,t=0}^d Q_{ri}Q_{sj}Q_{th}\left[uvw\atop rst\right]$.
┼ёыш ярЁрьхЄЁ ╩Ёхщэр $q_{ij}^h=0$, Єю $S_{ijh}(u,v,w)=0$.

╟рЇшъёшЁєхь тхЁ°шэ√ $u,v,w$ фшёЄрэЎшюээю Ёхуєы Ёэюую уЁрЇр $\Gamma$ фшрьхЄЁр $3$ ш яюыюцшь
$\{ijh\}=\left\{uvw\atop ijh\right\}$, $[ijh]=\left[uvw\atop ijh\right]$, $[ijh]'=\left[uwv\atop ihj\right]$,
$[ijh]^*=\left[vuw\atop jih\right]$ ш $[ijh]^\sim=\left[wvu\atop hji\right]$.  ┬ ёыєўр ї
$d(u,v)=d(u,w)=d(v,w)=2$ шыш $d(u,v)=d(u,w)=d(v,w)=3$ т√ўшёыхэшх ярЁрьхЄЁют $[ijh]'=\left[uwv\atop ihj\right]$,
$[ijh]^*=\left[vuw\atop jih\right]$ ш $[ijh]^\sim=\left[wvu\atop hji\right]$ (ёшььхЄЁшчрЎш  ьрёёштр ЄЁющэ√ї
ўшёхы яхЁхёхўхэшщ) ьюцхЄ фрЄ№ эют√х ёююЄэю°хэш , яючтюы ■∙шх фюърчрЄ№ эхёє∙хёЄтютрэшх уЁрЇр.

\section*{╤тющёЄтр уЁрЇр $\Gamma_3(u)$}

┬ ¤Єюь Ёрчфхых $\Gamma$  ты хЄё  фшёЄрэЎшюээю Ёхуєы Ёэ√ь уЁрЇюь ё ьрёёштюь яхЁхёхўхэшщ $\{55,54,2;1,1,54\}$. ╥юуфр
$\Gamma$ шьххЄ ёяхъЄЁ $55^{1}, 7^{1617}, -1^{110}, -8^{1408}$, $1+55+2970+110=3136$ тхЁ°шэ ш фєры№эє■ ьрЄЁшЎє ёюсёЄтхээ√ї
чэрўхэшщ
$$
Q=\left(\begin{array}{rrrr}
1 & 1617 & 110 & 1408 \\
1 & \frac{1029}{5} & -2 & -\frac{1024}{5} \\
1 & -\frac{49}{15} & -2 & \frac{64}{15} \\
1 & -\frac{147}{5} & 54 & -\frac{128}{5}
\end{array}\right).
$$

─рыхх, уЁрЇ $\Gamma_3$  ты хЄё  $56\times 56$-Ёх°хЄъющ ш юъЁхёЄэюёЄ№ тхЁ°шэ√ т $\Gamma_3$  ты хЄё  юс·хфшэхэшхь фтєї
шчюышЁютрээ√ї 55-ъышъ.

\begin{lem}\label{intnumbers} ╫шёыр яхЁхёхўхэшщ уЁрЇр $\Gamma$ Ёртэ√

$(1)$ $p^1_{11}=0$, $p^1_{12}=54$, $p^1_{22}=2808$, $p^1_{23}=108$, $p^1_{33}=2$;

$(2)$ $p^2_{11}=1$, $p^2_{12}=52$, $p^2_{13}=2$, $p^2_{22}=2811$, $p^2_{23}=106$, $p^2_{33}=2$;

$(3)$ $p^3_{12}=54$, $p^3_{13}=1$, $p^3_{22}=2862$, $p^3_{23}=54$, $p^3_{33}=54$.
\end{lem}

\proof ╧Ё ь√х т√ўшёыхэш , ёюуырёэю \cite[ыхььр 4.1.7]{BCN}.

╟рЇшъёшЁєхь тхЁ°шэ√ $u,v,w$ уЁрЇр $\Gamma$ ш яюыюцшь $\{ijh\}=\left\{uvw\atop ijh\right\}$,
$[ijh]=\left[uvw\atop ijh\right]$. ╧єёЄ№ $\Delta=\Gamma_2(u)$ ш $\Lambda=\Delta_2$. ╥юуфр $\Lambda$ --- Ёхуєы Ёэ√щ
уЁрЇ ёЄхяхэш $p^2_{22}=2811$ эр $k_2=2970$ тхЁ°шэрї.

\begin{lem}\label{211} ╧єёЄ№ $d(u,v)=2,d(u,w)=d(v,w)=1$. ╥юуфр ЄЁющэ√х ўшёыр яхЁхёхўхэшщ Ёртэ√:

$(1)$ $[122]=52$, $[132]=2$;

$(2)$ $[212]=52$, $[221]=53$, $[222]=r_{1}+2650$, $[223]=-r_{1}+108$, $[232]=-r_{1}+106$, $[233]=r_{1}$;

$(3)$ $[312]=2$, $[322]=-r_{1}+106$, $[323]=[332]=r_{1}$, $[333]=-r_{1}+2$,

уфх $r_{1}\in \{0, 1, 2\}$.
\end{lem}

\proof
╙яЁю∙хэш  ЇюЁьєы $(+)$.

\begin{lem}\label{223} ╧єёЄ№ $d(u,v)=d(u,w)=2,d(v,w)=3$.
╥юуфр ЄЁющэ√х ўшёыр яхЁхёхўхэшщ Ёртэ√:

$(1)$ $[112]=[121]=-r_{6}+1$, $[113]=[131]=r_{6}$, $[122]=r_{7}$, $[123]=[132]=r_{6}-r_{7}+51$, $[133]=-2r_{6}+r_{7}-49$;

$(2)$ $[212]=[221]=52$, $[213]=[231]=0$, $[222]=r_{5}+r_{6}+2756-r_{7}=r_{9}$,
$[223]=[232]=-2r_{6}+r_{7}+3$, $[233]=2r_{6}-r_{7}+102$;

$(3)$ $[312]=[321]=r_{6}+1$, $[313]=[331]=1-r_{6}$,  $[322]=5-2r_{6}$,
$[323]=[332]=r_{6}$, $[333]=1$,

уфх $r_{6}\in \{0,1\}$, $r_{7}\in \{49,...,52\}$.

─рыхх, ышсю $[222]=2706$ ш $[222]=2r_{6}+50=r_{7}$, ышсю $[222]=2707$, $2r_{6}+51=r_{7}$
ш $r_{6}=0$.
\end{lem}

\proof
╤ яюью∙№■ ъюья№■ЄхЁэ√ї єяЁю∙хэшщ ЇюЁьєы $(+)$ яюыєўшь ЁртхэёЄтр

$[112]=-r_{6}+1$, $[113]=r_{6}$, $[121]=-r_{4}-r_{8}+54$, $[122]=r_{7}$, $[123]=r_{4}-r_{7}+r_{8}-2$,
$[131]=r_{4}+r_{8}-53$, $[132]=r_{6}-r_{7}+51$, $[133]=-r_{4}-r_{6}+r_{7}-r_{8}+4$;

$[212]=r_{5}+r_{6}-r_{7}-r_{9}+2808$, $[213]=-r_{5}-r_{6}+r_{7}+r_{9}-2756$, $[221]=r_{4}$, $[222]=r_{9}$,
$[223]=-r_{4}-r_{9}+2811$, $[231]=-r_{4}+52$, $[232]=-r_{5}-r_{6}+r_{7}+3$, $[233]=r_{4}+r_{5}+r_{6}-r_{7}+50$;

$[312]=-r_{5}+r_{7}+r_{9}-2755$, $[313]=r_{5}-r_{7}-r_{9}+2757$, $[321]=r_{8}$, $[322]=-r_{7}-r_{9}+2861$,
$[323]=r_{7}-r_{8}+r_{9}-2755$, $[331]=-r_{8}+2$, $[332]=r_{5}$, $[333]=-r_{5}+r_{8}$,

уфх $r_{4}\in \{51, 52\}$, $r_{5}\in \{0,1,2\}$, $r_{6}\in \{0,1\}$, $r_{7}\in \{49,...,52\}$, $r_{8}\in \{0,1,2\}$, $r_{9}\in \{2705,...,2710\}$.

╬Єё■фр $2705\le [222]=r_{9}\le 2710$.

┬ $56\times 56$-Ёх°хЄъх $\Gamma_3$ яюфуЁрЇ $\Gamma_3(u)\cap \Gamma_3(v)$  ты хЄё  2-ъюъышъющ, яхЁхёхър■∙хщ
$\Gamma_3(w)$ яю хфшэёЄтхээющ тхЁ°шэх, яю¤Єюьє $[333]=-r_{5}+r_{8}=1$ ш $r_{5}\in \{0,1\}$.
└эрыюушўэю, $\Gamma_3(v)\cap \Gamma_3(w)$  ты хЄё  54-ъышъющ, яхЁхёхър■∙хщ
$\Gamma(u)$ эх сюыхх ўхь яю юфэющ тхЁ°шэх, яю¤Єюьє $[133]=-r_{4}-r_{6}+r_{7}-r_{5}+3\le 1$.

╤шььхЄЁшчрЎш . ─ы  ЄЁющъш тхЁ°шэ $(u,v,w)$ шьххь $[122]=r_{7}=r_{7}'$, $-r_{6}+1=[112]=[121]'=-r_{4}'-r_{5}'+53$
ш $r_{4}+r_{5}=r_{6}'+52$. ─рыхх, $[222]=r_{9}=r_{9}'$, $r_{6}=[323]=[332]'=r_{5}'$ ш $r_{4}=52$.
╬Єё■фр $[133]=-49-r_{6}+r_{7}-r_{5}\le 1$ ш $r_{7}\le r_{6}+r_{5}+50$.
╚ьххь $[313]=r_{5}-r_{7}-r_{9}+2757$, яю¤Єюьє $r_{7}+r_{9}\le r_{5}+2757$ ш $[213]=-r_{5}-r_{6}+r_{7}+r_{9}-2756\le 1-r_{6}$.

─юяєёЄшь, ўЄю $[213]=1$. ╥юуфр $r_{6}=0$, $[233]=r_{5}-r_{7}+102$  ш $r_{5}=r_{5}'=0$, яю¤Єюьє $r_{7}+r_{9}=2757$. ╥ръ ъръ
$[133]=r_{7}-r_{8}-48=r_{7}'-r_{8}'-48$, Єю $r_{8}=r_{8}'$. ─рыхх, $1=[112]=[121]'=r_{8}'$ ш $r_{8}=1$. └эрыюушўэю,
$[233]=r_{4}-r_{7}+50$, ёыхфютрЄхы№эю, $r_{4}=r_{4}'$. ╧ЁюЄштюЁхўшх ё Єхь, ўЄю $51=[212]=[221]'=r_{4}'=52$.

╟эрўшЄ, $[213]=0$ ш $r_{7}+r_{9}=r_{5}+r_{6}+2756$. ┼ёыш $r_{6}=0$, Єю $[233]=r_{5}-r_{7}+102$  ш $r_{5}=r_{5}'=1$.
┼ёыш цх $r_{6}=1$, Єю $[233]=r_{5}-r_{7}+103$  ш $r_{5}=r_{5}'=1$. ╧ю¤Єюьє т ы■сюь ёыєўрх $r_{5}=r_{6}$ ш яюыєўрхЄё 
ёшььхЄЁшчютрээ√щ ьрёёшт:

$[112]=[121]=-r_{6}+1$, $[113]=[131]=r_{6}$, $[122]=r_{7}$, $[123]=[132]=r_{6}-r_{7}+51$, $[133]=-2r_{6}+r_{7}-49$;

$[212]=[221]=52$, $[213]=[231]=0$, $[222]=r_{5}+r_{6}+2756-r_{7}=r_{9}$,
$[223]=[232]=-2r_{6}+r_{7}+3$, $[233]=2r_{6}-r_{7}+102$;

$[312]=[321]=r_{6}+1$, $[313]=[331]=1-r_{6}$,  $[322]=5-2r_{6}$,
$[323]=[332]=r_{6}$, $[333]=1$,

уфх $r_{6}\in \{0,1\}$, $r_{7}\in \{49,...,52\}$.

┼ёыш $[222]\ge 2708$, Єю $r_{7}\le 48+2r_{6}$ ш ЁртхэёЄтю $[133]=-2r_{6}+r_{7}-49$ фрхЄ яЁюЄштюЁхўшх.
┼ёыш $[222]=2705$, Єю $r_{7}=51+2r_{6}$, яЁюЄштюЁхўшх ё Єхь, ўЄю $[133]=-2r_{6}+r_{7}-49\le 1$.
╚Єръ, ышсю $[222]=2706$ ш $[222]=r_{5}+r_{6}+50=r_{7}$, ышсю $[222]=2707$, $r_{5}+r_{6}+51=r_{7}$
ш $r_{5}+r_{6}=0$.

\begin{lem}\label{221} ╧єёЄ№ $d(u,v)=d(u,w)=2,d(v,w)=1$.
╥юуфр ЄЁющэ√х ўшёыр яхЁхёхўхэшщ Ёртэ√:

$(1)$ $[122]=51$, $[123]=[132]=0$, $[133]=2$;

$(2)$ $[212]=[221]=51$, $[222]=2654$, $[223]=[232]=106$, $[233]=0$;

$(3)$ $[312]=[321]=2$, $[322]=102$, $[323]=[332]=2$, $[333]=0$.

─рыхх, $[222]=2654$.
\end{lem}

\proof
╤ яюью∙№■ єяЁю∙хэшщ ЇюЁьєы $(+)$ яюыєўшь ЁртхэёЄтр

$[122]=r_{3}$, $[123]=[132]=-r_{3}+51$, $[133]=r_{3}-49$;

$[212]=[221]=51$, $[222]=-r_{2}-r_{3}+2705$, $[223]=[232]=r_{2}+r_{3}+55$, $[233]=-r_{2}-r_{3}+51$;

$[312]=[321]=2$, $[322]=r_{2}+102$, $[323]=[332]=-r_{2}+2$, $[333]=r_{2}$,

уфх $r_{2}\in \{0, 1, 2\}$, $r_{3}\in \{49,50,51\}$.

┬ $56\times 56$-Ёх°хЄъх $\Gamma_3$ яюфуЁрЇ $\Gamma_3(u)\cap \Gamma_3(w)$  ты хЄё  2-ъюъышъющ, эх яхЁхёхър■∙хщ
$\Gamma_3(v)$. ╧ю¤Єюьє $[333]=r_{2}=0$.

╧юфёўшЄрхь ўшёыю $d$ ярЁ тхЁ°шэ $(y,z)$ эр ЁрёёЄю эшш $3$ т уЁрЇх $\Lambda$, уфх $y\in \left\{uv\atop 21\right\}$ ш
$z\in \left\{uv\atop 23\right\}$. ╤ юфэющ ёЄюЁюэ√, яю ыхььх 2 шьххь $[213]=0$ ш $d=0$.
╤ фЁєующ ёЄюЁюэ√, шьххь $d=52[233]=0$, яю¤Єюьє $[233]=0$ ш $r_{3}=51$. ╬Єё■фр ёыхфє■Є ЁртхэёЄтр шч
чръы■ўхэш  ыхьь√.

\begin{lem}\label{222} ╧єёЄ№ $d(u,v)=d(u,w)=(v,w)=2$.
╥юуфр ЄЁющэ√х ўшёыр яхЁхёхўхэшщ Ёртэ√:

$[111]=-r_{16}-r_{11}+1$, $[112]=r_{11}$, $[113]=r_{16}$, $[121]=r_{14}$, $[122]=-r_{12}-r_{14}+52$,
$[123]=r_{12}$, $[131]=-r_{14}+r_{16}+r_{11}$, $[132]=r_{12}+r_{14}-r_{11}$, $[133]=-r_{12}-r_{16}+2$;

$[211]=r_{15}+r_{16}-r_{17}+50$, $[212]=r_{17}$, $[213]=-r_{15}-r_{16}+2$, $[221]=-r_{15}-r_{16}+r_{17}-r_{10}+2$,
$[222]=r_{12}+r_{16}-r_{17}+r_{10}+2704$, $[223]=-r_{12}+r_{15}+104$, $[231]=r_{10}$,
$[232]=-r_{12}-r_{16}-r_{10}+106$, $[233]=r_{12}+r_{16}$;

$[311]=-r_{15}+r_{17}+r_{11}-50$, $[312]=-r_{17}-r_{11}+52$, $[313]=r_{15}$,
$[321]=-r_{14}+r_{15}+r_{16}-r_{17}+r_{10}+50$, $[322]=r_{14}-r_{16}+r_{17}-r_{10}+54$,
$[323]=-r_{15}+2$, $[331]=r_{14}-r_{16}-r_{10}-r_{11}+2$, $[332]=-r_{14}+r_{16}+r_{10}+r_{11}$, $[333]=0$,

уфх $r_{10},r_{12},r_{15}\in \{0,1,2\}$, $r_{11},r_{14},r_{16}\in \{0,1\}$, $r_{17}\in \{49,...,52\}$.

─рыхх, $2652\le [222]=r_{12}+r_{16}-r_{17}+r_{10}+2704\le 2659$.
\end{lem}

\proof
╤ яюью∙№■ ъюья№■ЄхЁэ√ї єяЁю∙хэшщ ЇюЁьєы $(+)$ яюыєўшь ЁртхэёЄтр:

$[111]=-r_{16}-r_{11}+1$, $[112]=r_{11}$, $[113]=r_{16}$, $[121]=r_{14}$, $[122]=-r_{12}-r_{14}+52$,
$[123]=r_{12}$, $[131]=-r_{14}+r_{16}+r_{11}$, $[132]=r_{12}+r_{14}-r_{11}$, $[133]=-r_{12}-r_{16}+2$;

$[211]=r_{15}+r_{16}-r_{17}+50$, $[212]=r_{17}$, $[213]=-r_{15}-r_{16}+2$, $[221]=-r_{15}-r_{16}+r_{17}-r_{10}+2$,
$[222]=r_{12}-r_{13}+r_{16}-r_{17}+r_{10}+2704$, $[223]=-r_{12}+r_{13}+r_{15}+104$, $[231]=r_{10}$,
$[232]=-r_{12}+r_{13}-r_{16}-r_{10}+106$, $[233]=r_{12}-r_{13}+r_{16}$;

$[311]=-r_{15}+r_{17}+r_{11}-50$, $[312]=-r_{17}-r_{11}+52$, $[313]=r_{15}$,
$[321]=-r_{14}+r_{15}+r_{16}-r_{17}+r_{10}+50$, $[322]=r_{13}+r_{14}-r_{16}+r_{17}-r_{10}+54$,
$[323]=-r_{13}-r_{15}+2$, $[331]=r_{14}-r_{16}-r_{10}-r_{11}+2$, $[332]=-r_{13}-r_{14}+r_{16}+r_{10}+r_{11}$, $[333]=r_{13}$,

уфх $r_{10},r_{12},r_{13},r_{15}\in \{0,1,2\}$, $r_{11},r_{14},r_{16}\in \{0,1\}$, $r_{17}\in \{49,...,52\}$.

┬ $56\times 56$-Ёх°хЄъх $\Gamma_3$ яюфуЁрЇ $\Gamma_3(u)\cap \Gamma_3(w)$  ты хЄё  2-ъюъышъющ, эх яхЁхёхър■∙хщ
$\Gamma_3(v)$. ╧ю¤Єюьє $[333]=r_{13}=0$. ╬Єё■фр ёыхфє■Є ЄЁхсєхь√х ЁртхэёЄтр.

╟рьхЄшь, ўЄю $r_{12}+r_{16}\le 2$, яю¤Єюьє $2652\le [222]=r_{12}+r_{16}-r_{17}+r_{10}+2704\le 2659$.
╦хььр фюърчрэр.
\medskip

─ы  ўшёыр ЁхсхЁ $e$ ьхцфє $\Lambda(v)$ ш $\Lambda_2(v)$ т уЁрЇх $\Lambda$ т√яюыэ ■Єё  эхЁртхэёЄтр
$424844=52\cdot 2654+106\cdot 2706\le e\le 52\cdot 2654+106\cdot 2707=424950$, яю¤Єюьє
$151.136\le 2810-\lambda\le 151.174$ ш $2658.826\le \lambda\le 2658.864$, уфх $\lambda$ --- ёЁхфэхх чэрўхэшх ярЁрьхЄЁр
$\lambda(\Lambda)$.

╟рьхЄшь, ўЄю ёЁхфэхх чэрўхэшх $\lambda(\Lambda)$ юўхэ№ сышчъю ъ тхЁїэхщ уЁрэшЎх 2659.

\section*{─юърчрЄхы№ёЄтю ЄхюЁхь√}\label{s4}

╧єёЄ№ $\Gamma$  ты хЄё  фшёЄрэЎшюээю Ёхуєы Ёэ√ь уЁрЇюь ё ьрёёштюь яхЁхёхўхэшщ
$\{55,54,2;1,1,54\}$. ┬ ¤Єюь Ёрчфхых ь√ фюърцхь ЄхюЁхьє \ref{Jgraph}.

╟рЇшъёшЁєхь тхЁ°шэ√ $u,v,w$ уЁрЇр $\Gamma$ ш яюыюцшь $\{ijh\}=\left\{uvw\atop ijh\right\}$,
$[ijh]=\left[uvw\atop ijh\right]$. ╧єёЄ№ $\Delta=\Gamma_2(u)$ ш $\Lambda=\Delta_2$. ╥юуфр $\Lambda$ --- Ёхуєы Ёэ√щ
уЁрЇ ёЄхяхэш $p^2_{22}=2811$ эр $k_2=2970$ тхЁ°шэрї.

\begin{lem}\label{sym}
╧єёЄ№ $d(u,v)=d(u,w)=(v,w)=2$. ╥юуфр т√яюыэ ■Єё  ёыхфє■∙шх єЄтхЁцфхэш :

$(1)$ тхЁэ√ ЁртхэёЄтр: $r_{11}=r_{14}'$, $r_{10}'=[213]=r_{12}^*$, $r_{15}+r_{16}+r_{10}'=2$ ш
$-r_{17}+52=r_{12}^*+r_{14}^*$;

$(2)$ хёыш $r_{10}'=2$, Єю $r_{10}=2$,

ышсю $r_{11}=0$ ш

$[111]=1$, $[112]=[113]=[121]=[131]=0$, $[122]=50$,
$[123]=[132]=2$, $[133]=0$;

$[211]=0$, $[212]=r_{17}=50$, $[213]=[231]=2$, $[221]=50$,
$[222]=2658$, $[223]=[232]=102$, $[233]=2$;

$[311]=0$, $[312]=[321]=2$, $[313]=0$, $[322]=102$,
$[323]=[332]=2$, $[331]=[333]=0$,

ышсю $r_{11}=1$ ш

$[111]=0$, $[112]=[121]=1$, $[113]=[131]=0$, $[122]=49$,
$[123]=[132]=2$, $[133]=0$;

$[211]=1$, $[212]=[221]=r_{17}=49$, $[213]=[231]=2$,
$[222]=2659$, $[223]=[232]=102$, $[233]=2$;

$[311]=0$, $[312]=[321]=2$, $[313]=[331]=0$,
$[322]=102$, $[323]=[332]=2$, $[333]=0$.
\end{lem}

\proof
╤шььхЄЁшчрЎш . ╥ръ ъръ $[111]=-r_{16}-r_{11}+1$, Єю $r_{16}+r_{11}=r_{16}'+r_{11}'=r_{16}^*+r_{11}^*=r_{16}^\sim+r_{11}^\sim$.
─рыхх, $[112]=r_{11}=r_{11}^*$, $[113]=r_{16}=r_{16}^*$, $[121]=r_{14}=r_{14}^\sim$, $[212]=r_{17}=r_{17}^\sim$,
$[313]=r_{15}=r_{15}^\sim$, $r_{11}=[112]=[121]'=r_{14}'$, $r_{10}'=[213]=r_{12}^*$.

╚ьххь $[122]=-r_{12}-r_{14}+52$, яю¤Єюьє $r_{12}+r_{14}=r_{12}'+r_{14}'$, $[233]=r_{12}+r_{16}=r_{12}'+r_{16}'$,
$r_{10}=[231]=[213]'=-r_{15}'-r_{16}'+2$ ш $r_{10}'+r_{15}+r_{16}=2$

└эрыюушўэю, $-r_{17}-r_{11}+52=[312]=[132]^*=r_{12}^*+r_{14}^*-r_{11}^*$ ш $-r_{17}+52=r_{12}^*+r_{14}^*$.
╙ЄтхЁцфхэшх (1) фюърчрэю.

─юяєёЄшь, ўЄю $r_{10}'=2$. ╥юуфр $r_{12}^*=2$, $r_{15}+r_{16}=0$ ш $r_{17}+r_{14}^*=50$. ╥ръ ъръ $[131]=-r_{14}+r_{11}$,
Єю т  ёыєўрх $r_{11}=0$ шьххь $r_{14}=0$, $[132]=r_{12}+r_{14}-r_{11}=r_{12}=r_{12}'$, $[211]=-r_{17}+50$, $[311]=r_{17}-50$.
╥хяхЁ№ $-r_{12}-r_{10}+106=[232]=[223]'=-r_{12}'+104=-r_{12}+104$, яю¤Єюьє $r_{10}=2$, $[212]=[221]=r_{17}=50$,
$r_{17}=[122]=-r_{12}+52$, $r_{12}=2$ ш

$[111]=1$, $[112]=[113]=[121]=[131]=0$, $[122]=50$,
$[123]=[132]=2$, $[133]=0$;

$[211]=0$, $[212]=r_{17}=50$, $[213]=[231]=2$, $[221]=50$,
$[222]=2658$, $[223]=[232]=102$, $[233]=2$;

$[311]=0$, $[312]=[321]=2$, $[313]=0$, $[322]=102$,
$[323]=[332]=2$, $[331]=[333]=0$.

┬ ёыєўрх $r_{11}=1$ шьххь $[233]=r_{12}=r_{12}'$, яю¤Єюьє $r_{12}=[123]=[132]=[132]=r_{12}+r_{14}-1$, $r_{14}=1$,
$-r_{17}+51=[312]=[321]'=-r_{17}'+r_{10}'+49$, $r_{17}=r_{17}'$, $r_{10}=2$ ш $-r_{17}+51=[312]=[132]^*=r_{12}^*$.
╥ръ ъръ $[222]=r_{12}+2657$, Єю $r_{12}=r_{12}^*$,  $r_{12}+r_{17}=51$ ш

$[111]=0$, $[112]=[121]=1$, $[113]=[131]=0$, $[122]=49$,
$[123]=[132]=2$, $[133]=0$;

$[211]=1$, $[212]=[221]=r_{17}=49$, $[213]=[231]=2$,
$[222]=2659$, $[223]=[232]=102$, $[233]=2$;

$[311]=0$, $[312]=[321]=2$, $[313]=[331]=0$,
$[322]=102$, $[323]=[332]=2$, $[333]=0$.

╥хяхЁ№ $[212]=[221]=r_{17}=r_{17}^*$, $r_{17}=[122]=-r_{12}+51$ ш $r_{12}=2$.
╦хььр фюърчрэр.
\medskip

╤ыхфє■∙р  шфх  ┬шфрыш яючтюы хЄ ёюъЁрЄшЄ№ фюърчрЄхы№ёЄтю ЄхюЁхь√.  ╧єёЄ№ $d(u,v)=2$.
╥юуфр $[u]$ ёюфхЁцшЄ яю юфэюьє ¤ыхьхэЄє т ърцфющ ёЄЁюъх ш ърцфюь ёЄюысЎх $55\times 55$-Ёх°хЄъш 
$\Gamma_3-(\{u\}\cup \Gamma_3(u))$. ╧юыюцшь $\{x,y\}=[u]\cap \Gamma_3(v)$. ╥юуфр шьххЄё  хфшэёЄтхээр  тхЁ°шэр
$w$ т яюфуЁрЇх $\Gamma_3(x)\cap \Gamma_3(y)-\{v\}$. ╧ЁюЄштюЁхўшх ё Єхь, ўЄю $[133] = 2$ (ёь. ыхььє \ref{221}) фы  
ы■сющ ЄЁющъш тхЁ°шэ $(u,v,w)$ ё $d(u,v)=d(u,w)=2,d(v,w)=1$.  

╥хюЁхьр \ref{Jgraph}  фюърчрэр.

\bigskip

╤тхфхэш  юс ртЄюЁрї:

\smallskip

╠└╒═┼┬ └ыхъёрэфЁ └ыхъёххтшў,

\smallskip

╨╬╤╤╚▀, у. ┼ърЄхЁшэсєЁу, єы. ╤.╩ютрыхтёъющ, 16, 620990,

╚эёЄшЄєЄ ьрЄхьрЄшъш ш ьхїрэшъш шь. ═.═. ╩Ёрёютёъюую ╙Ё╬ ╨└═,

╨╬╤╤╚▀, у. ┼ърЄхЁшэсєЁу, єы. ╠шЁр, 19, 620002,

╙Ёры№ёъшщ ЇхфхЁры№э√щ єэштхЁёшЄхЄ,

х-mail: makhnev@imm.uran.ru

\end{document}